\documentclass[pdftex,letter,11pt]{article}

\usepackage{url}
\makeatletter
\def\url@leostyle{%
  \@ifundefined{selectfont}{\def\UrlFont{\sf}}{\def\UrlFont{\small\ttfamily}}}
\makeatother
\usepackage[english]{babel}
\usepackage{color}
\usepackage{natbib}
\bibpunct{[}{]}{,}{n}{,}{,}
\usepackage{indentfirst}
\usepackage{fancyhdr}
\usepackage[dvips]{graphicx}
\usepackage[figuresright]{rotating}
\usepackage{amsthm}
\usepackage[toc,page]{appendix}% for help see http://www.tex.ac.uk/cgi-bin/texfaq2html?label=appendix

\usepackage{epsfig} % for postscript graphics files
\usepackage{mathptmx} % assumes new font selection scheme installed
\usepackage{times} % assumes new font selection scheme installed
\usepackage{latexsym}
\usepackage{amssymb}
\usepackage[fleqn,tbtags]{amsmath} % for help read amsldoc.pdf
\usepackage{amsfonts}
\usepackage{makeidx}
\usepackage{subfigure}
\usepackage{rotating}%The rotating package by Sebastian Rahtz and Leonor Barroca tries to make the interface for rotation somewhat simpler. It defines
\usepackage[all]{xy}
\usepackage{ae} %per vocali accentate; added by Giuseppe Bilotta: Type1 Fonts with T1 enc
\swapnumbers
\newtheorem{theorem}{Theorem}[section]
\newtheorem{lemma}[theorem]{Lemma}
\newtheorem{corollary}[theorem]{Corollary}
\newtheorem{proposition}[theorem]{Proposition}
\newtheorem{definition}[theorem]{Definition}

\newtheorem{conjecture}[theorem]{Conjecture}

\def\N{{\mathbb{N}}}
\def\R{{\mathbb{R}}}

\def\C{{\mathbb{C}}}
\def\tr{{\mathrm{tr\,}}}
\def\F{{\mathcal F}}

\def\P{{\mathcal P}}

\def\PF{{\mathcal P}({\cal F})}

\def\rhoF{\rho({\cal F})}
\def\nxn{n\! \times \!n\,}
\def\2x2{2\! \times \!2\,}
\def\DLIF{\textrm{DLI}(\F)}

\def\Cn{\C^{\, n}}

\title{A note on the Joint Spectral Radius}

\author{Antonio Cicone \thanks{Universit\'a degli Studi dell'Aquila,
      via Vetoio 1, 67100, L'Aquila, Italy\newline {E-mail: \tt\small antonio.cicone@univaq.it}}}

\begin{document}

\maketitle
\date{}

Last two decades have been characterized by an increasing interest
in the analysis of the maximal growth rate of long products
generated by matrices belonging to a specific set/family. The
maximal growth rate can be evaluated considering a generalization of
the spectral radius of a single matrix to the case of a set of
matrices.

This generalization can be formulated in many different ways,
nevertheless in the commonly studied cases of bounded or finite
families all the possible generalizations coincide in a unique value
that is usually called \emph{joint spectral radius} or simply
\emph{spectral radius}. The joint spectral radius, however, can
prove to be hard to compute and can lead even to undecidable
problems. We present in this paper all the possible
generalizations of the spectral radius, their properties and the
associated theoretical challenges.

From an historical point of view the first two generalizations of
spectral radius, the so--called joint and common spectral radius,
were introduced by Rota and Strang in the three pages paper
``\emph{A note on the joint spectral radius}'' published in $1960$
\cite{RotStr}. After that more than thirty years had to pass before
a second paper was issued on this topic: in $1992$ Daubechies and
Lagarias \cite{DauLag1} published ``\emph{Sets of matrices all
infinite products of which converge}'' introducing the generalized
spectral radius, conjecturing it was equal to the joint spectral
radius (this was proven immediately after by Berger and Wang
\cite{BerWan}) and presenting examples of applications. From then on
there has been a rapidly increasing interest on this subject and the
more years pass the more the number of mathematical branches and
applications directly involved in the study of these quantities
increases \cite{Blondel}.

The study of infinite products convergence properties proves to be
of primary interest in a variety of contexts:

\noindent Nonhomogeneous Markov chains, deterministic construction
of functions and curves with self-similarities under changes in
scale like the von Koch snowflake and the de Rham curves, two-scale
refinement equations that arise in the construction of wavelets of
compact support and in the dyadic interpolation schemes of
Deslauriers and Dubuc  \cite{DauLag1,Pro4}, the asymptotic behavior
of the solutions of linear difference equations with variable
coefficients \cite{GugZen00,GugZen01,GugZen03}, coordination of
autonomous agents \cite{JadLinMor03,Moreau,GhoshLee}, hybrid systems
with applications that range from intelligent traffic systems to
industrial process control \cite{BlTs99}, the stability analysis of
dynamical systems of autonomous differential equations
\cite{BraTon}, computer--aided geometric design in constructing
parametrized curves and surfaces by subdivision or refinement
algorithms \cite{MiPra89,CavaMi}, the stability of asynchronous
processes in control theory \cite{Tsi87}, the stability of
desynchronised systems \cite{Koz90}, the analysis of magnetic
recording systems and in particular the study of the capacity of
codes submitted to forbidden differences constraints
\cite{MOrSi01,BJP06}, probabilistic automata \cite{Pat}, the
distribution of random power series and the asymptotic behavior of
the Euler partition function \cite{Pro4}, the logarithm of the joint
spectral radius appears also in the context of discrete linear
inclusions as the Lyapunov indicator \cite{Bara88,Gurv}.
%- Lagarias Wang 1995 and also arise in studying the
%dynamical complementarity problem in the theory of stochastic
%networks [see Kozyakin et al. (1993)].
For a more extensive and detailed list of applications we refer the
reader to the Gilbert Strang's paper ``\emph{The Joint Spectral
Radius}'' \cite{Str} and to the doctoral theses by Jungers and Theys
\cite{Ju09,Theys}.

The paper develops as following: in Section \ref{sec:Notation} we
give notation and terminology used throughout this paper; Section \ref{sec:FrameworkJSR} presents first a case
of study associated with the asymptotic behavior analysis of the
solutions of linear difference equations with variable coefficients,
further, it contains the definitions and properties of all the
possible generalizations of spectral radius for a set of matrices,
in particular the irreducibility, nondefectivity and finiteness properties are discussed.%;
%last Section is an overview on the so--called joint spectral
%subradius, Lyapunov exponent and p--radii and their relation to the
%joint spectral radius.

\section{Terminology, notation and basic properties}\label{sec:Notation}

In this Section we provide notation, terminology, definitions and
properties which are employed in this paper.

We use the expression $\N_{0}$ meaning the set of natural numbers,
included zero. All the matrices and vectors that we consider have real or complex
entries. We denote the \index{conjugate transpose}\emph{conjugate
transpose} of an $m$--by--$n$ matrix $X=[x_{ij}]$ by
$X^*=[\bar{x}_{ji}]$, while the simple \index{transpose}
\emph{transpose} as $X^T=[x_{ji}]$.

For $p\in [1,\infty)$, the
$l^p$~norm of a vector $w\in \mathbb{C}^n$ is given by
$
\|{w}\|_{p}= \sqrt[p]{\sum_{i=1}^n |w[i]|^p}
$

In particular:
\begin{description}
  \item[$l^1$ -- The sum norm]\index{vector norm!sum
norm}\index{vector norm!$l^1$|see{vector norm, sum norm}}
  $
   \|w\|_1=\sum_i|w[i]|
  $
  \item[$l^2$ -- The Euclidean norm]\index{vector norm!Euclidean norm}\index{vector norm!$l^2$|see{vector norm, Euclidean norm}}
  $
   \|w\|_2=\sqrt{\sum_{i=1}^n |w[i]|^2}=\sqrt{w^*w}
  $
\item[$l^\infty$ -- The max norm]\index{vector norm!max norm}\index{vector
norm!$l^\infty$|see{vector norm, max norm}}
$
\|{w}\|_{\infty}= \max_{j=1,\ldots,n} |w[j]|.
$
\end{description}
%For a square matrix $A\in \mathbb{C}^{\nxn}$ and for $p\in
%[1,\infty]$, $\|A\|_{p}$ is the associated induced norm given by
%$
% \|A\|_{p}=\max_{\|x\|_{p}=1}\|Ax\|_{p}
%$

If $A$ is a square matrix, its \index{characteristic polynomial}
\emph{characteristic polynomial} is $p_{A}(t):=\det(tI-A)$, where
$\det$ stands for determinant \cite[Section 0.3]{horn}; the
(complex) zeroes of $p_{A}(t)$ are the
\index{eigenvalue}\emph{eigenvalues} of $A$. A complex number
$\lambda$ is an eigenvalue of $A$ if and only if there are nonzero
vectors $x$ and $y$ such that $Ax = \lambda x$ and $y^{\ast}A =
\lambda y^{\ast}$; $x$ is said to be an \emph{eigenvector} (more
specifically, a \index{eigenvector!right}\emph{right eigenvector})
of $A$ associated with $\lambda$ and $y$ is said to be a
\index{eigenvector!left}\emph{left eigenvector} of $A$ associated
with $\lambda$. The set of all the eigenvalues of $A$ is called the
\index{spectrum of a matrix}\emph{spectrum} of $A$ and is denoted by
$\sigma(A)$. The \emph{determinant}\index{determinant} of $A$, $\det
A$, is equivalent to the product of all its eigenvalues. If the
spectrum of $A$ does not contain $0$ the matrix is said
\emph{nonsingular}\index{nonsingular matrix} ($A$ nonsingular if and
only if $\det A\neq 0$). The \index{spectral radius}\emph{spectral
radius} of $A$ is the nonnegative real number $\rho(A) = \max\left\{
|\lambda|: \lambda \in \sigma(A)\right\}$. If $\lambda\in
\sigma(A)$, its \index{multiplicity!algebraic}\emph{algebraic
multiplicity} is its multiplicity as a zero of $p_{A}(t)$; its
\index{multiplicity!geometric}\emph{geometric multiplicity} is the
maximum number of linearly independent eigenvectors associated with
it. The geometric multiplicity of an eigenvalue is never greater
than its algebraic multiplicity. An eigenvalue whose algebraic
multiplicity is one is said to be
\index{eigenvalue!simple}\emph{simple}. An eigenvalue $\lambda$ of
$A$ is  said to be \emph{semisimple} \index{eigenvalue!semisimple}
if and only if rank$(A-\lambda I) = $rank$(A-\lambda I)^2$ i.e.
$\lambda$ has the same geometric and algebraic multiplicity. If the
geometric multiplicity and the algebraic multiplicity are equal for
every eigenvalue, $A$ is said to be \index{nondefective
matrix}\emph{nondefective}, otherwise is
\emph{defective}\index{defective|see{nondefective
matrix}}\label{defective1}.

We let $e_{1}$ indicate the first column of the identity matrix $I$:
$e_{1}=[1~0~\cdots~0]^T$. We let $e =[1~1~\cdots~1]^T$ denote the
all--ones vector. Whenever it is useful to indicate that an identity
or zero matrix has a specific size, e.g., $r$--by--$r$, we write
$I_r$ or $0_r$.

Two vectors $x$ and $y$ of the same size are \emph{orthogonal}
\index{vectors!orthogonal} if $x^{\ast}y = 0$. The \emph{orthogonal
complement} \index{orthogonal complement} of a given set of vectors
is the set (actually, a vector space) of all vectors that are
orthogonal to every vector in the given set.

An $n$--by--$r$ matrix $X$ has \index{vectors!orthonormal}
\emph{orthonormal} columns if $X^*X = I_r$. A square matrix $U$ is
\index{matrix!unitary} \emph{unitary} if it has orthonormal columns,
that is, if $U^{\ast}$ is the inverse of $U$.

A square matrix $A$ is a \emph{projection} \index{projection} if
$A^2=A$.

A square matrix $A$ is \index{matrix!row--stochastic}
\emph{row--stochastic} if it has real nonnegative entries and
$Ae=e$, which means that the sum of the entries in each row is 1;
$A$ is \index{matrix!column--stochastic} \emph{column--stochastic}
if $A^T$ is row--stochastic. We say that $A$ is
\index{matrix!stochastic} \emph{stochastic} if it is either
row--stochastic or column--stochastic.\label{def:stoch}

The \emph{direct sum} \index{direct sum} of $k$ given square
matrices $X_1, \ldots, X_k$ is the block diagonal matrix
\[
\left[
\begin{array}{ccc}
X_1  & \cdots  & 0  \\
\vdots & \ddots  & \vdots \\
0 & \cdots & X_k%
\end{array}
\right] = X_1 \oplus \cdots \oplus X_k.
\]

The $k$--by--$k$ \emph{Jordan block} \index{Jordan block} with
eigenvalue $\lambda$ is
\[
J_{k}(\lambda)=\left[
\begin{array}
[c]{cccc}%
\lambda & 1 &  & 0\\
& \ddots & \ddots & \\
&  & \ddots & 1\\
&  &  & \lambda
\end{array}
\right]  \text{,} \quad J_1(\lambda)=[\lambda] \text{.}
\]

\noindent Each square complex matrix $A$ is similar to a direct sum
of Jordan blocks, which is unique up to permutation of the blocks;
this direct sum is the \emph{Jordan canonical form} \index{Jordan
canonical form} \label{def:JordCanForm} of $A$. The algebraic
multiplicity of $\lambda$ as an eigenvalue of $J_k(\lambda)$ is $k$;
its geometric multiplicity is $1$. If $\lambda$ is a semisimple
eigenvalue of $A$ with multiplicity $m$, then the Jordan canonical
form of $A$ is $\lambda I_{m} \oplus J$, in which $J$ is a direct
sum of Jordan blocks with eigenvalues different from $\lambda$; if
$\lambda$ is a simple eigenvalue, then $m=1$ and the Jordan
canonical form of $A$ is $[\lambda] \oplus J$. $A$ is
diagonalizable, i.e. its Jordan canonical form is given by a
diagonal matrix, if and only if is nondefective.

In a block matrix, the symbol $\bigstar$ denotes a block whose
entries are not required to take particular values.

We consider $A^0=I$. A matrix $B$ is said to be \emph{normal} \index{matrix!normal} if
$B\,B^*=B^*\!B$, \emph{unitary} \index{matrix!unitary} if
$B\,B^*=B^*\!B=I$, \index{matrix!Hermitian} \emph{Hermitian} if
$B=B^*$. Hermitian and unitary matrices are, by definition, normal
matrices.

A \emph{proper subset}\index{subset!proper} of a set $\cal A$ is a
set $\cal B$ that is strictly contained in $\cal A$. This is written
as $\cal B \varsubsetneq \cal A$.

Besides the Jordan canonical form, we need to introduce an additional matrix factorization, the so--called
\emph{singular value decomposition} (in short svd): Given a square
matrix $A\in \mathbb{C}^{\nxn}$ with rank $k\leq n$, there always
exists a diagonal matrix $\Lambda\in \mathbb{R}^{\nxn}$ with
nonnegative diagonal entries $\sigma_1\geq \sigma_2\geq \cdots \geq
\sigma_k > \sigma_{k+1} = \cdots = \sigma_n = 0$ and two unitary
matrices $U,V \in \mathbb{C}^{\nxn}$ such that $A =U\Lambda V^*$,
which is defined as the singular value decomposition of
$A$\index{singular value decomposition}. The matrix $\Lambda =
\textrm{diag}(\sigma_1,\ldots, \sigma_n)$ is always uniquely
determined and $\sigma_1^2\geq \cdots\geq \sigma_n^2$ correspond to
the eigenvalues of the Hermitian matrix $AA^*$. Values ${\sigma_1,
\ldots, \sigma_n}$ are \index{singular values}the so--called
\emph{singular values} of $A$.

The \index{trace}\emph{trace} of an $\nxn$--matrix $A$, denoted by
$\tr(A)$, is given by the sum of the diagonal elements of $A$,
$\tr(A) = \sum_{i=0}^n a_{ii}$, and it is also equal to the sum of
all the eigenvalues in the spectrum of $A$, $\tr (A) =
\sum\limits_{\lambda\in\sigma(A)}\lambda$.

The \emph{spectral radius}  \index{spectral radius} of a square
matrix $A\in \mathbb{C}^{\nxn}$ is defined as
\begin{equation}
\rho(A) = \max\left\{ |\lambda|:\lambda \in \sigma(A)\right\}
\end{equation}
It is easy to prove that $\rho(A^k) = (\rho(A))^k$ for every
$k\in\N$ and, thus, given a generic power $k$ of the matrix $A$, the
value $(\rho(A^k))^{1/k}$ is just equal to the spectral radius of
the matrix.

It is possible to characterize the spectral radius using the
trace\index{spectral radius!via trace} of the matrix. Since
$\lambda^k\in\sigma(A^k)$ for every eigenvalue $\lambda\in\sigma(A)$
and for every $k\in\N$, it follows that
\begin{equation}
|\tr(A^k)|^{1/k}=\rho(A)\left|\sum\limits_{\lambda\in\sigma(A)}
\lambda^k/(\rho(A))^k\right|^{1/k}
\end{equation}
which converges to $\rho(A)$ as $k\rightarrow \infty$
\begin{equation}
\rho(A)=\lim\limits_{k\rightarrow
\infty}\left|\tr(A^k)\right|^{1/k}\!\!\!\!\!\!\!\!\!.
\end{equation}

For a square matrix $A$ and for $p\in [1,\infty]$, $\|A\|_{p}$ is
the \index{matrix norm!induced} matrix norm induced by the
corresponding $p$--vector norm. The
induced matrix norms are sometimes defined as \emph{operator
norms}\index{operator norm} \cite[Definition 5.6.3]{horn}. Among the
induced matrix norms we will make use of the following
\begin{description}
  \item[The maximum column--sum norm]\index{matrix norm!maximum column--sum
norm}
  $\|A\|_1=\displaystyle{\max_{\|x\|_{1}=1}}\|Ax\|_{1}=\displaystyle{\max_{j}\sum_i|a_{i,j}|}$
  \item[The spectral norm]\index{matrix norm!spectral norm}
  $\|A\|_2=\displaystyle{\max_{\|x\|_{2}=1}}\|Ax\|_{2}=\sigma_1(A)=\sqrt{\rho(AA^*)}$
  \item[The maximum row--sum norm]\index{matrix norm!maximum row--sum
norm}
  $\|A\|_\infty=\displaystyle{\max_{\|x\|_{\infty}=1}}\|Ax\|_{\infty}=\displaystyle{\max_{i}\sum_j}|a_{i,j}|$
\end{description}
Every induced matrix norm $\|\cdot \|_*$ is submultiplicative i.e.
$\|AB\|_*\leq \|A\|_*\|B\|_*$ for every square matrix $A$ and $B$.

Another family of induced matrix norms are the\index{matrix
norm!ellipsoidal}\index{vector norm!ellipsoidal}\label{ellipsNorm}
\emph{ellipsoidal norms}. Let us consider an Hermitian
\index{positive definite|see{matrix, positive definite}}
\index{matrix!positive definite} positive definite matrix $P\succ 0$
(i.e. $P$ is a nonsingular Hermitian matrix such that $x^*Px > 0$
for all nonzero $x \in \mathbb{C}^n$ or, equivalently, $P$ is a
Hermitian matrix such that all its eigenvalues are strictly
positive). The \emph{vector ellipsoidal norm} is defined as
\begin{equation}
\|x\|_P=\sqrt{x^*Px}.
\end{equation}

The corresponding induced matrix norm is given by
\begin{equation}\label{eq:EllNorm}
\|A\|_P=\max_{\|x\|_{P}=1}\|Ax\|_{P}=\max_{\sqrt{x^*Px}=1}\sqrt{x^*A^*PAx}
\end{equation}

Recalling that \cite[Corollary 7.2.9]{horn} $P$ is positive definite
if and only if there exists a nonsingular upper triangular matrix
$T\in \mathbb{C}^{\nxn}$, with strictly positive diagonal entries,
such that $P = T^*T$, which is defined as the \index{Cholesky
decomposition}\emph{Cholesky decomposition} of  $P$, we can rewrite
$
\|A\|_P = \max_{\sqrt{x^*T^*Tx}=1}\sqrt{x^*A^*T^*TAx} =
\max_{\|Tx\|_2=1}\|TAx\|_2
$
\noindent and if we rename $y=Tx$, $T$ by construction is
nonsingular so $x=T^{-1}y$, we get
\begin{equation}\label{ellnorm}
\|A\|_P=\max_{\|y\|_2=
1}\|TAT^{-1}y\|_2=\|TAT^{-1}\|_2=\sqrt{\rho(TAT^{-1}(TAT^{-1})^*)}
\end{equation}
Since $T$ is nonsingular and remembering that the spectrum of a
matrix is invariant under\index{similarity transformation}
\emph{similarity transformation}, two matrices $M$ and $T^{-1}MT$
have the same eigenvalues, counting multiplicity. So from
(\ref{ellnorm}) we obtain that
\begin{equation}\label{ellnorm2}
\|A\|_P=\sqrt{\rho(TAP^{-1}A^*T^*)}=\sqrt{\rho(AP^{-1}A^*P)}
\end{equation}

Given a generic power $k$ of the matrix $A$, the value
$\|A^k\|^{1/k}$ is defined as the \emph{normalized
norm}\index{norm!normalized} of the matrix, in the sense that is
normalized with respect to the length of the
product.\index{normalized norm|see{norm, normalized}}

Given the family $\F=\{A_i\}_{i\in {\cal I}}$ of complex square
$\nxn$--matrices, with ${\cal I}$ a set of indices, $\F$ is defined
\index{bounded|see{set, bounded}} \index{set!bounded} \emph{bounded}
\label{def:bounded} if it does exist a constant $C < +\infty$ such
that $\sup\limits_{i\in {\cal I}}\|A_i\|\leq C $. While we define
the set \index{finite|see{set, finite}} \index{set!finite}
\emph{finite} if it is constituted by a finite number of matrices.
Trivially a finite set is always bounded.

A matrix $A$ is said to be \index{nondefective
matrix}\emph{nondefective}\label{defective2} if and only if its
Jordan canonical form is diagonal i.e. each eigenvalue of $A$ is
semisimple or, equivalently, it has geometric multiplicity equal to
algebraic multiplicity, otherwise $A$ is defined \emph{defective}.
In this paper we deal with a weaker condition of
nondefectivity: a matrix $A$ is said to be
\index{nondefective!weakly}\emph{weakly nondefective} if and only if
the eigenvalues of $A$ with modulus equal to the spectral radius,
i.e. with maximal modulus, are semisimple, if it is not the case the
matrix is defined \index{weakly defective|see{nondefective,
weakly}}\emph{weakly defective}. Using the Jordan canonical form of
$A$ it is easy to prove that, whenever $\rho(A)> 0$, defined
$A^\ast=A/\rho(A)$, $A$ is weakly nondefective if and only if powers
$(A^\ast)^k$ are bounded for every $k \geq 1$.

\label{def:defective}From now on, for the sake of simplicity and to
be coherent with the literature on the spectral radius of sets, we
use the expressions \emph{strongly nondefective} and \emph{strongly
defective} in place of nondefective and defective, whereas we make
use of the words \emph{nondefective} and \emph{defective} meaning
weakly nondefective and weakly defective.

Let us now recall basic relations between spectral radius and matrix
norms:
\begin{theorem}[{\cite[Theorem 5.6.9]{horn}}]\label{th:normDomin}
If $\|\cdot\|$ is any matrix norm on $\mathbb{C}^{\nxn}$ and if
$A\in \mathbb{C}^{\nxn}$, then
$
\rho(A)\leq \|A\|.
$
\end{theorem}

Furthermore
\begin{lemma}[{\cite[Lemma 5.6.10]{horn}}]\label{InfNorm}
Let $A\in \mathbb{C}^{\nxn}$, for every $\epsilon> 0$ there is a
matrix norm $\|\cdot\|_\epsilon$ such that
\begin{equation}\label{eq:infNorm}
\rho(A)\leq \|A\|_\epsilon\leq \rho(A)+\epsilon
\end{equation}
\end{lemma}

The spectral radius of $A$ is not itself a matrix or vector norm,
but if we let $\epsilon\rightarrow 0$ in (\ref{eq:infNorm}) we have
that $\rho(A)$ is the greatest lower bound for the values of all
matrix norms of $A$
\begin{equation}\label{eq:RhoasInfNorm}
\rho(A)=\inf_{\|\cdot\|\in {\cal N}}\|A\|
\end{equation}
where ${\cal N}$ denotes the set of all possible induced matrix
norms\index{spectral radius!via $\inf$ on norms} (the so--called
operator norms).

Spectral radius allows to characterize \index{matrix!convergent}
\emph{convergent matrices}, i.e. those matrices whose successive
powers tends to zero:
\begin{theorem}[{\cite[Theorem 5.6.12]{horn}}]\label{convergenceCondition}
    Let $A\in\mathbb{C}^{\nxn}$, then
    $
    \lim\limits_{k\rightarrow \infty}A^k=0\qquad\Leftrightarrow\qquad \rho(A)<1
    $
\end{theorem}
As a Corollary of the previous Theorem we have the so--called
\index{Gelfand's formula}\emph{Gelfand's formula}:
\begin{corollary}[{\cite[Corollary 5.6.14]{horn}}]\label{GelfandFormula}
Let $\|\cdot\|$ be any matrix norm on $\mathbb{C}^{\nxn}$, then
\begin{equation}\label{eq:GelfandFormula}
\rho(A)= \lim\limits_{k\rightarrow \infty}\|A^k\|^{1/k} \ \ \textrm{
for all } \ \ A\in \mathbb{C}^{\nxn}
\end{equation}
\end{corollary}
\noindent The Gelfand's formula gives us two information:
\begin{itemize}
  \item the spectral radius of $A$ represent the asymptotic growth rate of the
normalized norm of $A^k$: $\|A^k\|^{1/k}\sim \rho(A) \textrm{ as
}k\rightarrow \infty$
  \item the normalized norm $\|A^k\|^{1/k}$
  \index{spectral radius!via normalized norm} can
  be used to approximate the spectral
  radius and in the limit for
  $k\rightarrow\infty$ the two quantities coincide.
\end{itemize}

Given a row--stochastic matrix $A$ its maximum row--sum matrix norm
is equal $1$ by definition of row--stochasticity. By Theorem
\ref{th:normDomin}, choosing as matrix norm the maximum row--sum, we
have that for every stochastic matrix $A$
\begin{equation}\label{eq:SRstoch}
\rho(A)\leq\|A\|_\infty=1
\end{equation}
The row--stochasticity of $A$ can be formulated also as
\begin{equation}\label{eq:MatStoch}
A e =e
\end{equation}
with $e$ the all--ones vector and $\lambda = 1$ the eigenvalue of
$A$ associated with the right eigenvector $x=e$. So we have that
$\rho(A)=1$. Remembering that the
set of stochastic matrices is closed under matrix multiplication, we
observe that the very same result can be proved also using the
Gelfand formula: choosing as matrix norm the maximum row--sum we
have that $\|A^k\|_\infty^{1/k}=1$ for every integer $k\geq 1$.

In the following we generalize all these notions to the case of a
family of matrices.

For a systematic discussion of a broad range of
matrix analysis issues, see \cite{horn}.

\section{Framework}\label{sec:FrameworkJSR}
\subsection{A case of study}\label{subsec:exampleJSR}

Given a stable discrete time system we want to analyze its
robustness with respect to perturbations not a priori quantifiable.

Let us consider the system
\begin{equation}\label{eq:stabSys}
x(k+1) = A_{0}\,x(k), \qquad k \in \N_{0}.
\end{equation}
with  $x(0) \in \Cn$ and $A_{0} \in \C^{\,\nxn}$ such that the
system is asymptotically stable, i.e. $\rho(A_0) < 1$ (ref Theorem
\ref{convergenceCondition}). We consider the perturbed system given
by time--varying perturbations
\begin{equation}\label{eq:pertSys}
x(k+1) = \left( A_{0} + \sum\limits_{i=1}^{p} \delta_{i}(k)\,A_{i}
\right)\,x(k), \qquad k \in \N_{0}.
\end{equation}
The matrices $\{ A_i \}_{i=1}^{p}$ are known, but the perturbations
$\{ \delta_i (k) \}_{i=1}^{p}$ are not. The perturbations  may
depend on incomplete modeling, neglect of dynamics or measurement
uncertainty. We are interested to know if a stability result for the
theoretical model (\ref{eq:stabSys}) holds also for the real system
(\ref{eq:pertSys}).

The perturbed system (\ref{eq:pertSys}) can be regarded as a first
order system of difference equations with variable coefficients
\begin{equation}\label{eq:sd}
x(k+1) = Y_{i_k}\,x(k), \qquad k \in \N_{0}.
\end{equation}
where $x(0) \in \Cn$ and $Y_{i_k} \in \C^{\,\nxn}$ is an element of
the following family
\begin{equation}\label{eq:Family}
\F_{\alpha} = \left\{ A_{0} + \sum\limits_{i=1}^{p}
\delta_{i}\,A_{i} \ \bigg| \ \| \delta \| \le \alpha \right\}
\end{equation}
where $\delta = \left( \delta_1 \ \delta_2 \ \cdots \ \delta_p
\right)^T$ and the bound on the uncertainties is known. This kind of
problems arise in several contexts such as when applying numerical
methods to non--autonomous systems of differential equations.

From a point of view of robustness or worst case analysis the goal
is to determine the largest uncertainty level $\alpha^{\ast}$ such
that for every $\alpha < \alpha^{\ast}$ the system remains stable
(see e.g. \cite{Wirt02}).

If the sequence of matrices $Y_{i_k}$ is known, for $k \ge 0$, then
the solution of (\ref{eq:sd}) is given by
\begin{equation}\label{eq:SysSol}
x(k+1) = P_{k}\,x(0), \qquad \mbox{with} \ \ P_{k} =
{\prod\limits_{j=1}^{k} Y_{i_j}}, \ k \ge 1
\end{equation}
where asymptotic stability may be studied directly (although this is
not an easy task in general). Nevertheless we want to study the case
where the sequence $\{ Y_{i_k} \}_{k \ge 1}$ is not known a priori
and may be whatever.

\begin{definition} [Uniform asymptotic stability -- u.a.s.]\label{def:uas}
We say that {\rm (\ref{eq:sd})} \index{u.a.s.|see{uniform asymptotic
stability}} is {\em uniformly asymptotically stable}\index{uniform
asymptotic stability} if
\begin{equation}\label{eq:uas}
\lim\limits_{k \rightarrow \infty} x(k) = 0
\end{equation}
for any initial $x(0)$ and any sequence $\{ Y_{i_k} \}_{k \ge 1}$ of
elements in $\F_{\alpha}$.
\end{definition}

It is easy to prove that Definition \ref{def:uas} is equivalent of
requiring that any possible left product $Y_{i_k}\cdot
Y_{i_{k-1}}\!\!\cdot\ldots\cdot Y_{i_1}$ of matrices from
$\F_{\alpha}$ vanishes as $k \rightarrow \infty$.

We observe that in the context of the discrete linear inclusions
some authors refer to the uniform asymptotic stability as
\emph{absolute asymptotic stability}\index{absolute asymptotic
stability} \cite{Gurv,Theys}.

For the single matrix case we have that u.a.s. holds if and only if
the spectral radius of the matrix is strictly less than one, while
for the general case of a family of matrices $\F$ we are driven to
the problem of computing the {\em joint spectral radius} of $\F$.
The intrinsic difficulty in exploiting this quantity is due to the
non--commutativity of matrix multiplication.

\subsection{Definitions and properties}\label{subsec:Def&PropJSR}
\subsubsection{Definitions}\label{subsubsec:Def}

From now on we consider always  complex square $n\! \times
\!n$--matrices and submultiplicative norms if not differently
specified. Let $\F=\left\{A_i\right\}_{i\in {\cal I}}$ be a family
of matrices, ${\cal I}$ being a set of indices.

\noindent For each $k=1,2,\ldots $, consider the set ${\mathcal
P}_k(\F)$ of all possible products of length $k$ whose factors are
elements of $\F$, that is $\P_k(\F)=\{A_{i_1}\cdot\ldots \cdot
A_{i_k}\ | \ i_1,\ldots,i_k \in {\cal I}\}$ and set
\begin{equation}\label{eq:Msgr}
\PF = \bigcup\limits_{k \ge 1} \P_k(\F)
\end{equation}
to be the \index{semigroup!multiplicative}\emph{multiplicative
semigroup} associated with $\F$. While, defined ${\mathcal
P}_0(\F):=I$, we have
\begin{equation}\label{eq:Mmon}
\P^*(\F) = \bigcup\limits_{k \ge 0} \P_k(\F)
\end{equation}
the \index{monoid!multiplicative}\emph{multiplicative monoid}
associated with $\F$.

We present four different generalizations of the concept of spectral
radius of a single matrix to the case of a family of matrices $\F$.

The first generalization is due to Rota and Strang, in the seminal
paper \cite{RotStr} published in $1960$ they presented the
generalization of the notion of spectral radius as limit of the
normalized norm of a single matrix\index{spectral radius!via
normalized norm}:
\begin{definition}[Joint Spectral Radius -- jsr]\label{def:JSR} If
$\|\cdot\|$ is any matrix norm on $\mathbb{C}^{\nxn}$, consider
$\widehat\rho_k(\F):=\sup_{P\in \P_k(\F)} \|P\|^{1/k},
\qquad k\in \N$ i.e. the supremum among the normalized norms of all
products in $\P_k(\F)$, and define the \index{spectral
radius!joint}{\em joint spectral radius} of $\F$ as
\begin{equation}\label{eq:JSR}
\textrm{jsr}\,(\F)=\widehat\rho(\F)=\lim_{k\rightarrow \infty} \
\widehat\rho_k(\F)
\end{equation}
\end{definition}
The joint spectral radius does not depend on the matrix norm chosen
thanks to the equivalence between matrix norms in finite dimensional
spaces.

We observe that in the discrete linear inclusions literature the
logarithm of the joint spectral radius is sometimes called Lyapunov
indicator \cite{Bara88}.

In $1992$ Daubechies and Lagarias \cite{DauLag1} introduced the
generalized spectral radius as a generalization of the $\limsup$
over all the spectral radii $\rho(A^k)^{1/k}$, $k\geq 1$, which are,
trivially, always equal to $\rho(A)$.
\begin{definition}[Generalized Spectral Radius -- gsr]\label{def:GSR}
Let $\overline\rho(\cdot)$ denote the spectral radius of an $n\!
\times\! n $--matrix, consider
$\overline\rho_k(\F):=\sup_{P\in \P_k(\F)} \rho(P)^{1/k},
\qquad k\in \N$ i.e. the supremum among the spectral radii of all
products in $\P_k(\F)$ normalized taking their $k$--th root, and
define the \index{spectral radius!generalized}{\em generalized
spectral radius} of $\F$ as
\begin{equation}\label{eq:GSR}
\textrm{gsr}\,(\F)=\overline\rho(\F)=\limsup\limits_{k\rightarrow
\infty} \ \overline\rho_k(\F)
\end{equation}
\end{definition}

For this two definitions it has been proved by Daubechies and
Lagarias \cite{DauLag1,DauLag2} the following
\begin{proposition}[Four members inequality]\label{prop:fourMember}
For any set of matrices $\F$ and any $k\geq 1$
\begin{equation}\label{eq:fourMember}
\overline\rho_k(\F)\leq \overline\rho(\F)=\textrm{gsr}\,(\F)\leq
\textrm{jsr}\,(\F)=\widehat\rho(\F) \leq \widehat\rho_k(\F)
\end{equation}
independently of the submultiplicative norm used to define
$\widehat\rho_k(\F)$.
\end{proposition}
As a consequence of this we have that:
\begin{equation}\label{eq:JSRinf}
\widehat\rho(\F)=\inf\limits_{k \geq 1}\, \widehat\rho_k(\F)
\end{equation}
\begin{equation}\label{eq:GSRsup}
\overline\rho(\F)=\sup\limits_{k \geq 1}\, \overline\rho_k(\F)
\end{equation}
For the first equality see also \cite[Lemma 1.2]{Ju09}; for the
second one,  since $\rho(M^k)=\rho(M)^k$ for every $k\in\N$ and
considering that by  definition of $\limsup$
\begin{equation*}
\limsup\limits_{k\rightarrow \infty} \
\overline\rho_k(\F)=\inf_{k\geq 1} \ \sup_{n\,\geq k} \
\overline\rho_n(\F),
\end{equation*}
\noindent if it does exist a finite product $P\in\P_r(\F)$,
$r\in\N$, such that $\rho(P)^{1/r}=\overline\rho(\F)$, then, for
every $m\in\N$, $\rho(P^{\,m})^{1/mr}=\overline\rho(\F)$ and, thus,
$\sup_{n\,\geq k}\,\overline\rho_n(\F)=\overline\rho(\F)$ for every
$k\in\N$. This last equality is valid also if it does not exists
such a finite product, in fact in this case the $\sup$ is achieved
only for $n\rightarrow \infty$. So in both cases it results
$\inf_{\,k\geq 1} \, \sup_{n\,\geq k} \
\overline\rho_n(\F)=\sup_{k\geq 1} \ \overline\rho_k(\F)$, i.e.
equation (\ref{eq:GSRsup}) holds true.

A third definition has been introduced by Chen and Zhou in 2000
\cite{CheZho} and is based on a generalization of the formula
associating the spectral radius of a matrix with its
trace\index{spectral radius!via trace}:
\begin{definition}[Mutual Spectral Radius -- msr]\label{def:MSR}
Let $\tr (P)$ be the trace of the product $P\in \P_k(\F)$ then
$\sup\limits_{P\in\P_k(\F)}|\tr (P)|$ is the maximal absolute value
among all the traces of the products of length $k$. Define the
\index{spectral radius!mutual}{\em mutual spectral radius} of ${\cal
F}$ as
\begin{equation}\label{eq:MSR}
\textrm{msr}\,(\F)=\limsup\limits_{k\rightarrow \infty}
\sup\limits_{P\in\P_k(\F)}|\tr (P)|^{1/k}
\end{equation}
\end{definition}

We present now the last characterization of the spectral radius of a
family of matrices. For bounded sets (ref Section
\ref{sec:Notation}) it is possible to generalize the concept,
express in equation (\ref{eq:RhoasInfNorm}), of spectral radius as
the $\inf$ over the set of all possible induced matrix norms of $A$.
\index{spectral radius!via $\inf$ on norms}
\begin{definition}[Common Spectral Radius -- csr]\label{def:CSR}
Given a norm $\|\cdot \|$ on the vector space $\Cn$ and the
corresponding induced matrix norm, we define
\begin{equation}\label{eq:normF}
\|\F\|:=\sup_{i\in {\cal I}}\|A_i\|
\end{equation}
where we assume that the family $\F$ is bounded. We define the
\index{spectral radius!common}{\em common spectral radius} of ${\cal
F}$ (see \cite{RotStr} and \cite{Els}) as
\begin{equation}\label{eq:CSR}
\textrm{csr}\,(\F)=\widetilde{\rho}(\F)=\inf_{\|\cdot\|\in {\cal
N}}\|\F\|
\end{equation}
where ${\cal N}$ denotes the set of all possible induced matrix
norms.
\end{definition}
This definition was first introduced by Rota and Strang in $1960$
\cite{RotStr} and re--introduced 35 years later by Elsner
\cite{Els}.

In the case of bounded sets, it is possible to prove that the four
characterizations we presented coincide.
\begin{theorem}[The Complete Spectral Radius Theorem]\label{th:prop2}
For a bounded family $\F$ the following equalities hold true
\begin{equation}\label{eq:defn}
\textrm{gsr}\,(\F)=\textrm{jsr}\,(\F)=\textrm{csr}\,({\cal
F})=\textrm{msr}\,(\F)
\end{equation}
\end{theorem}

The equality of gsr and jsr was conjectured by Daubechies and
Lagarias and it was proven by Berger and Wang \cite{BerWan}, Elsner
\cite{Els}, Chen and Zhou \cite{CheZho}, Shih et al. \cite{ShWuPa}.
For the equality of csr and jsr we refer the reader to the seminal
work of Rota and Strang \cite{RotStr} or again \cite{Els}. Chen and
Zhou \cite{CheZho} proved the last equality.

We observe that the first equality is the generalization of the
Gelfand's formula (Corollary \ref{GelfandFormula}) to the case of a
family of matrices.

Another observation is that even though the joint and generalized
spectral radius can be defined also for unbounded families the first
equality does not hold in general. Consider for example the
unbounded family:
$$
\F=\left\{\left(
                        \begin{array}{cc}
                          1 & 1 \\
                          0 & 1 \\
                        \end{array}
                      \right),\ldots,\left(
                        \begin{array}{cc}
                          1 & n \\
                          0 & 1 \\
                        \end{array}
                      \right),\ldots\right\}
$$
For this family since every product of the two matrices is upper
triangular with ones in the main diagonal it is evident that
$\overline\rho(\F)=1$ and obviously $\widehat\rho(\F)=+\infty$ since
the family is unbounded (see \cite{Theys} for details and
\cite{DauLag1} for another example).

We observe also that Gurvits in \cite{Gurv} give a counterexample to
the first equality in the case of two operators in an infinite
dimensional Hilbert space.

From now on and if not differently specified we will always consider
bounded sets of matrices. Theorem \ref{th:prop2} implies that we can
simply refer to the \index{spectral radius of a family}{\em spectral
radius} $\rho(\cal F)$ of the family of matrices $\F$.

\begin{definition}[Trajectory]\label{def:traj}
Given a family $\F= \left\{A_i \right\}_{i\in {\cal I}}$, we define,
for an arbitrary nonzero vector $x\in \Cn$, the \index{trajectory}
{\em trajectory}
\begin{equation}\label{eq:traj}
{\cal T}\left[\F,x\right]= \left\{P\,x \ \mid P \in \PF\right\}
\end{equation}
as the set of vectors obtained by applying all the products $P$ in
the multiplicative semigroup $\PF$ to the vector $x$.
\end{definition}

\begin{definition}[Discrete linear inclusion]\label{def:DLI}
The \index{discrete linear inclusion}\emph{discrete linear
inclusion} is the set of all the trajectories associated with all
the possible vectors in $\Cn$. This set is denoted by $\DLIF$.
\end{definition}

\subsubsection{Properties}\label{subsubsec:PropJSR}
We resume now properties valid for the spectral radius of a bounded
set of matrices $\F=\left\{A_i\right\}_{i\in {\cal I}}$
\begin{description}
  \item[1. Multiplication by a scalar:] For any set $\F$ and for any
number $\alpha\in\C$
\begin{equation}\label{eq:scalarInvariance}
\rho(\alpha \F) = |\alpha|\rho(\F)
\end{equation}

\item[2. Continuity:] The joint spectral radius is
continuous in its entries as established by Heil and Strang
\cite{HeiStra}. Wirth has proved \cite{Wirt02} that the joint
spectral radius is even locally Lipschitz continuous on the space of
compact irreducible sets of matrices (an explicit formula for the
related Lipschitz constant has been evaluated by Kozyakin
\cite{Koz10}).

\item[3. Powers of the family:] For any set $\F$ and for any $k \geq 1$
$$\rho(\F^k) \leq \rho^k(\F)$$

  \item[4. Invariance under similarity:] The spectral
radius of the family is invariant under similarity transformation,
so for any set of matrices $\F$, and any invertible matrix $T$
\begin{equation}\label{eq:invariance}
\rho(\F)=\rho(T\F T^{-1})
\end{equation}
This because to any product $A_1\cdot\ldots\ \cdot A_k \in {\mathcal
P}_k(\F)$ corresponds a product $T\cdot A_1\cdot\ldots\cdot A_k\cdot
T^{-1} \in \P_k(T\F T^{-1})$ with equal spectral radius.

  \item[5. Conjugate or transposed family:] The conjugate or transposed
  family (family obtained taking the conjugate/transpose of every
  matrix in the original set) has the same spectral radius as the original
  one \cite[Lemma 5.1]{GWZ05}
  \begin{equation}\label{eq:conj}
\rho(\F^*)=\rho(\F) \\
\rho(\F^T)=\rho(\F)
\end{equation}

\item[6. Block triangular matrices:]\label{blockTriang}
Given a family of block upper triangular matrices
\begin{equation*}
\F=\left\{\left(
\begin{array}{cc}
A_i & B_i \\
0 & C_i \\
\end{array}
\right)\right\}_{i\in {\cal I}}
\end{equation*}
we have that
\begin{equation}\label{eq:blockJSR}
\rho(\F) = \max \left\{\rho(\left\{A_i\right\}_{i\in {\cal I}}),
\rho(\left\{C_i\right\}_{i\in {\cal I}})\right\}.
\end{equation}
This is due to the closure, with respect to the multiplication, of
block upper triangularity \cite[Lemma II (c)]{BerWan}. Clearly the
same holds for lower triangular matrices. This result generalizes to
the case of more than two blocks on the diagonal.

\item[7. Closure and convex hull:]
The closure and the convex hull of a set have the same spectral
radius of the original set
\begin{equation}\label{eq:conv&clos}
\rho(\textrm{conv}\F)=\rho(\textrm{cl}\,\F)=\rho(\F)
\end{equation}
This result was first obtained by Barabanov in $1988$ \cite{Bara88}.
An alternative proof, given by Theys in \cite[page 17]{Theys}, is
based on the common spectral radius definition (\ref{eq:CSR}) and
the property
\begin{equation}\label{eq:eqSup}
    \sup_{A_i\in\F}\|A_i\|= \sup_{A_i\in\textrm{conv}\F}\|A_i\| =
    \sup_{A_i\in\textrm{cl}\,\F}\|A_i\|.
\end{equation}

\item[8. Uniform asymptotic stability characterization {\cite[Theorem I
(b)]{BerWan}}:]\label{item:u.a.s.Characterization} For any bounded
set of matrices $\F$ and for any $k\geq 1$, all matrix products
$P\in\P_k(\F)$ converge to the zero matrix as $k\rightarrow\infty$,
i.e. $\!\!\F$ is uniformly asymptotically stable (ref page
\pageref{def:uas}), if and only if $\rhoF < 1$.

In other words the spectral radius of the family of matrices $\F$
gives information about the uniform asymptotic stability of the
associated dynamical system $\DLIF$, defined on page
\pageref{def:DLI}.

\item[9. Product boundedness {\cite[Theorem I (a)]{BerWan}}:]
Given a bounded set of matrices $\F$, if products $P\in\P_k(\F)$,
$k\in\N$, converge as $k\rightarrow\infty$. Then, the multiplicative
monoid $\P^*(\F)$ defined in (\ref{eq:Mmon}) is bounded and $\rhoF
\leq 1$.

The opposite implication is not true in general:

Given a defective family with $\rhoF = 1$, products $P\in\P_k(\F)$,
$k\in\N$, explode for $k\rightarrow\infty$ by Definition
\ref{def:defectiveF}.

\noindent We return on this aspect in \cite{JSREC}.

\item[10. Special cases:]\label{specialCases}$\,$

\begin{enumerate}
  \item\label{item:stoc}
Recalling that the set of stochastic matrices is closed under matrix
multiplication and that every stochastic matrix has spectral radius
equal $1$ (ref Section \ref{sec:Notation}), if the matrices in $\F$
are all stochastic then the spectral radius of the family is exactly
$1$.

  \item\label{item:maxRho}
  If the matrices in $\F$ are all upper--triangular, if they can be
simultaneously upper--triangularized, if all the matrices in $\F$
commutes or, more in general, if the Lie algebra\index{Lie
algebra}\index{Lie algebra!solvable} associated with the set of
matrices is solvable (commutative families are \emph{Abelian Lie
algebras} which\index{Lie algebra!Abelian} are always solvable), if
they are all symmetric or, more in general, if they are all normal
or, finally, if they can be simultaneously normalized, then
\begin{equation}\label{eq:maxOverAll}
\rho(\F) = \max\limits_{A_i\in \F} \left\{\rho(A_i)\right\}
\end{equation}
For more details see \cite{Gurv,Gri96,CoHe,Theys,Ju09}.

  \item\label{item:c1}
  If $\F=\left\{A,A^*\right\}$ then $\rhoF=\rho(AA^*)^{1/2}=\sigma_1(A)$ i.e. the largest singular value of
  A. In fact \cite[Proposition 6.20]{Theys} using the four members inequality (\ref{eq:fourMember}) for $k=2$ we have
  \begin{equation}\label{eq:c1}
  \rho(AA^*)^{1/2}=\sigma_1(A)=\|AA^*\|_2^{1/2}
  \end{equation}

  \item({\cite{Moss} and \cite[Theorem 4]{GMV10}}).\label{item:Th4}
Consider the family $\F=\{A, B\}$ with
$$
A:=\left(
     \begin{array}{cc}
       a & b \\
       c & d \\
     \end{array}
   \right),\quad
   B:=\left(
     \begin{array}{cc}
       a & -b \\
       -c & d \\
     \end{array}
   \right)\quad
   a,b,c,d\in \R
$$
The joint spectral radius of the family $\F$ is given by
$$
\rho(\F) =\left\{ \begin{array}{l}
                          \rho(A)=\rho(B) \quad\textrm{if}\quad b\, c\geq0 \\
                          \sqrt{\rho(AB)} \quad\textrm{if}\quad b\, c<0 \\
                        \end{array}\right.
$$

\item\label{item:Th5} ({\cite{Moss} and \cite[Theorem 5]{GMV10}}).
Consider the family $\F=\{A, B\}$ with
$$
A:=\left(
     \begin{array}{cc}
       a & b \\
       c & d \\
     \end{array}
   \right),\quad
   B:=\left(
     \begin{array}{cc}
       d & c \\
       b & a \\
     \end{array}
   \right)\quad
   a,b,c,d\in \R
$$
The joint spectral radius of the family $\F$ is given by
$$
\rho(\F) =\left\{ \begin{array}{l}
                          \rho(A)=\rho(B) \quad\textrm{if}\quad |a-d|\geq |b-c| \\
                          \sqrt{\rho(AB)} \quad\textrm{if}\quad |a-d|< |b-c| \\
                        \end{array}\right.
$$

  \item\label{item:absF}
Let $|\F|$ be the family of matrices obtained from $\F$ as follows:
\[
A = [ a_{ij} ] \in \F \quad \longrightarrow \quad |A| =[ |a_{ij}| ]
\in |\F|.
\]
Then
\begin{equation}\label{eq:JSRabsF}
    \rho(|\F|)\geq\rho(\F)
\end{equation}

From the previous result and the four members inequality
(\ref{eq:fourMember}) we have that
\begin{equation*}
\overline\rho_k(\F)\leq\rho(\F)\leq\rho(|\F|)
\end{equation*}

So if $P \in \P_{k} (\F)$, $k\in\N$, is such that $\rho(P)^{1/k} =
\rho \left( | \F | \right)$, then
\begin{equation*}
\rho \left( \F \right) = \rho(P)^{1/k}.
\end{equation*}
\end{enumerate}

\item[11. Non--algebraicity:]\index{non--algebraicity}
Any set composed of $k$ real $\nxn$--matrices can be seen as a point
in the space $\R^{\,kn^2}$. Therefore, a subset of $\R^{\,kn^2}$ is
a set of $k$--tuples of $\nxn$--matrices. Given a subset of
$\R^{\,kn^2}$ this is defined \index{semi--algebraic set}
\emph{semi--algebraic} if it is a finite union of sets that can be
expressed by a finite list of polynomial equalities and
inequalities. Kozyakin \cite{Koz90} has shown that, for all $k,n\geq
2$, the set of points $x\in\R^{\,kn^2}$ such that $\rho(x)<1$ is not
semi--algebraic and, for all $k,n\geq 2$, the set of points
$x\in\R^{\,kn^2}$ corresponding to a bounded semigroup $\P(x)$ is
not semi--algebraic (the original paper by Kozyakin contains a flaw
and the correction has been published by the same author only in
Russian. For a corrected version in English we refer the reader to
the Doctoral work of Theys \cite[Section 4.2]{Theys}). In practice
in the general case, given a discrete linear inclusion $\DLIF$,
there is no procedure involving a finite number of operations that
allows to decide whether $\DLIF$ is uniformly asymptotically stable
or not i.e. the uniform asymptotic stability of $\DLIF$ is in
general hard to determine.

\item[12. NP--hardness:]\label{item:NP-hard}
In \cite{TsitBlon97} Tsitsiklis and Blondel proved that, given a set
of two matrices $\F$ and unless P $=$ NP, the spectral radius
$\rhoF$ is not polynomial--time approximable. This holds true even
if all nonzero entries of the two matrices are constrained to be
equal. Let us recall that the function $\rhoF$ is
\emph{polynomial--time approximable} \index{polynomial--time
approximable} if there exists an algorithm $\rho^*(\F,\epsilon)$,
which, for every rational number $\epsilon > 0$ and every set of
matrices $\F$ with $\rhoF > 0$, returns an approximation of $\rhoF$
with a relative error of at most $\epsilon$ (i.e. such that
$|\rho^{\ast}-\rho|\leq\epsilon\rho$) in time polynomial in the bit
size of $\F$ and $\epsilon$ (if  $\epsilon= p/q$, with $p$ and $q$
relatively prime numbers, its \index{bit size} \emph{bit size} is
equal to $\log(pq)$); however there are algorithms which are
polynomial either in the bit size of $\F$ or in $\epsilon$. We
conclude that the computation of the spectral radius of a set of
matrices is in general \index{NP--hard} \emph{NP--hard} and,
consequently, it is NP--hard to decide the stability of all products
of a set of matrices (for a survey of NP--hardness and
undecidability we refer the reader to \cite{BlTs00}). We observe
here that Gurvits in \cite{Gurv2} provides a polynomial--time
algorithm for the case of binary matrices.

\item[13. Undecidability:]\label{item:undec}
A decision problem is a problem which output is binary and can be
interpreted as ``yes'' or ``not''. For instance the problem of
deciding whether an integer matrix is nonsingular is a decision
problem. Since the nonsingularity\index{nonsingular matrix} can be
checked, for example, by computing the determinant of the matrix and
comparing it to zero it is a \emph{decidable
problem}\index{decidable problem}, i.e. a problem for which there
exists an algorithm that always halts with the right answer. But
there are also problems for which this kind of algorithm does not
exist, these are \emph{undecidable problems}.\index{undecidable
problem}

Given a set of matrices $\F$:
\begin{itemize}
  \item The problem of determining if the semigroup $\PF$ is bounded is
undecidable
  \item The problem of determining if $\rhoF\leq 1$ is undecidable
\end{itemize}
These two results, which remain true even if $\F$ contains only
rational entries \cite{BlTs00b,BloCan03}, teach us that does not
exist any algorithm allowing to compute the spectral radius of a
generic set $\F$ in finite time.

It is still unknown if it does exist in the generic case an
algorithm that, given a finite set of matrices $\F$, decides whether
$\rhoF< 1$. Such an algorithm would allow to check the uniform
asymptotic stability of the dynamical system ruled by the generic
set $\F$. In the following we discuss the relation between this kind
of algorithm and the so--called finiteness property.
\end{description}

The actual computation of $\rho(\F)$ is an important problem in
several applications, as we mentioned in the introduction of the
present paper. According to the previous properties of
non--algebraicity, NP--hardness and undecidability the problem
appears quite difficult in general.

However, this is not reason enough for declaring the problem
intractable and refraining from further research. As we discover in
the next subsection the existence of an s.m.p. for the family (i.e.
a product in the semigroup $\PF$ with particular properties) allows
in the general case to evaluate exactly the spectral radius of a
family making use of the Definition \ref{def:CSR} as an actual
computational tool. In order to do this we need the $\inf$ in
equation (\ref{eq:CSR}) to be a $\min$, but this is always true for
irreducible families.

\subsection{Irreducibility, nondefectivity and finiteness property}\label{subsec:Irr&Ndef}

When the $\inf$ in (\ref{eq:CSR}) is a $\min$ we say that the family
$\F=\left\{A_i\right\}_{i\in {\cal I}}$ admits an \emph{extremal
norm}.
\begin{definition}[Extremal norm]\label{def:extrNorm}
A norm $\|\cdot \|_*$ satisfying the condition
$$\rhoF=\|\F\|_*:=\sup\limits_{i\,\in {\cal I}}\|A_i\|_*$$
is said to be \index{norm!extremal}{\em extremal} for the family
$\F$ (for an extended discussion see \cite{Wirt05}).
\end{definition}

Equivalently a norm $\|\cdot \|_*$ is called extremal for a given
set $\F=\{A_i\}_{i\in {\cal I}}$ if it satisfies $\|A_i \|_*\leq
\rho(\F)$ for every $i\,\in {\cal I}$.

From Proposition \ref{prop:fourMember} it is clear that, for a given
norm, this inequality cannot be strict simultaneously for all the
matrices in the set.

Given a bounded family $\F=\{A_i\}_{i\in {\cal I}}$ of $n\! \times\!
n $--matrices with $\rho(\F)> 0$, the \emph{normalized family}
\index{normalized family}is given by
\begin{equation}\label{eq:NormFam}
{\cal F^{\ast}} = \left\{A_i/\rho(\F) \right\}_{i\in {\cal I}}
\end{equation}
with spectral radius $\rho(\F^{\ast})= 1$ and $\P
\left(\F^{\ast}\right)$ is the associated multiplicative semigroup
(ref equation (\ref{eq:Msgr}) on page \pageref{eq:Msgr}).

The definition of (weakly) defective matrix, given in section
\ref{sec:Notation}, extends to bounded families of matrices as
follows:
\begin{definition}[Defective and Nondefective Families]\label{def:defectiveF}
A bounded family $\F$ of $n\! \times\! n $--matrices is said to be
\index{defective family}\emph{defective} if the corresponding
normalized family ${\cal F^{\ast}}$ is such that the associated
semigroup $\P (\F^{\ast})$ is an unbounded set of matrices.
Otherwise, if either $\rhoF=0$ or $\rhoF>0$ with $\P
\left(\F^{\ast}\right)$ bounded, then the family $\F$ is said to be
\index{nondefective family}\emph{nondefective}.
\end{definition}

The following result can be found, for example, in \cite{RotStr} and
\cite{BerWan}:
\begin{proposition}\label{prop:nondef-extr}
A bounded family $\F$ of $n\! \times\! n $--matrices admits an
extremal norm $\|\cdot \|_*$ if and only if it is nondefective.
\end{proposition}

As previously mentioned we want to make use of Definition
\ref{def:CSR} as an actual computational tool for the spectral
radius $\rhoF$. To do this we need to ensure that the family admits
an \emph{extremal norm} i.e. we have to check the defectivity or
nondefectivity of the set $\F$.

Strictly connected to defectivity of a family there is the concept
of reducibility.
\begin{definition}[Reducible and Irreducible families]\label{def:decomp}
A bounded family $\F=\{A_i\}_{i\in {\cal I}}$ of $n\! \times\! n
$--matrices is said to be \index{reducible family}{\em reducible} if
there exist a nonsingular $n\! \times\! n $--matrix $M$ and two
integers $n_1,n_2\geq 1$, $n_1+n_2=n$, such that
\begin{equation}\label{eq:red}
M^{-1}A_iM=\left( \begin{array}{cc}
A_i^{(11)} & A_i^{(12)} \\
O & A_i^{(22)}
\end{array} \right) \quad\textrm{ for all }\quad i\in {\cal I}
\end{equation}
where the blocks $A_i^{(11)}$, $A_i^{(12)}$, $A_i^{(22)}$ are
$n_1\!\times n_1$--, $n_1\!\times n_2$-- and $n_2\!\times
n_2$--matrices, respectively. On the contrary, if a family ${\cal
F}$ is not reducible, then it is said to be {\em irreducible}.
\index{irreducible family}
\end{definition}
Irreducibility means that only the trivial subspaces ${0}$ and $\Cn$
are invariant under all the matrices of the family ${\cal F}$.
Otherwise $\F$ is called reducible. The concept of irreducibility
was introduced in the joint spectral radius theory by Barabanov in
\cite{Bara88}, where he named irreducible families \emph{nonsingular
sets}\index{nonsingular set|see{irreducible family}}.

We observe that some authors refer to reducibility as
decomposability\index{decomposable family|see{reducible family}}
(irreducibility as non--decomposability\index{non--decomposable
family|see{irreducible family}}) in order to avoid confusion with
the notion of reducibility commonly used in linear algebra \cite[Definition 6.2.21]{horn}.

An immediate consequence of irreducibility of $\F$ is that $\rho(\F)
> 0$, in fact, in this case the semigroup $\P(\F)$ is
irreducible and, therefore, does not consist of nilpotent elements,
by the Levitzky Theorem \cite{Lev}. So we can always normalize an
irreducible set of matrices $\F$ by $\rho(\F)$ obtaining a set with
generalized spectral radius equal to~$1$.

Another consequence of irreducibility of a family is stated in the
next Theorem and its Corollary, which follow easily from the
Barabanov's construction of extremal norms for irreducible families
of matrices \cite{Bara88}.

\begin{theorem}[{\cite[Lemma 4]{Els}}]\label{th:def-red}
    If a bounded family $\F$ of $n\times n$--matrices is
    defective, then is reducible.
\end{theorem}
and, therefore,
\begin{corollary}[{\cite{Bara88,Els,Pro3}}]
If a bounded family of matrices is irreducible then it is
nondefective, i.e. it does exist an extremal norm for the family.
\end{corollary}

In Figure \ref{fig:Mspace} it is represented the space ${\cal B}$ of
bounded families of matrices in $\C^{\,\nxn}$. This space can be
split into the set of the reducible families ${\cal R}$ and its
complement ${\cal I_R}$, the set of the irreducible ones.
%Let us consider three
%characterization of families of matrices: strongly nondefective,
%nondefective and defective. Trivially strongly nondefectivity
%implies nondefectivity by definition (ref page
%\pageref{def:defective}), but the opposite is not always true so,
%dealing with complementary sets, we have that the class of the
%defective families ${\cal D}$ is a proper subset of the set of the
%strongly defective families ${\cal S_D}$. Furthermore Theorem
%\ref{th:def-red} and the observation that there exist reducible
%matrices that are nondefective implies that the class ${\cal D}$ is
%a proper subset of ${\cal R}$ i.e. ${\cal D} \varsubsetneq {\cal
%R}$.
Families of matrices can be nondefective or defective: the set
${\cal D}$ of the defective families is a proper subset of ${\cal
R}$ i.e. ${\cal D} \varsubsetneq {\cal R}$. In fact Theorem
\ref{th:def-red} implies that a defective family is always
reducible, but the opposite implication is not necessarily true. For
example, for $n\geq 2$ all single families $\F=\{A\}$ are clearly
reducible as the Jordan canonical form proves, but not necessarily
defective. The set of nondefective families ${\cal N_D}$, the
complement of ${\cal D}$ in ${\cal B}$, is denoted by grey dots.
% or strongly defective.

About the dimension of set ${\cal D}$ and ${\cal R}$ Maesumi
\cite{Mae3} proposed the following conjecture
\begin{conjecture}
Reducible (decomposable) matrix sets form a set of measure zero in
the corresponding space of matrices. Defective matrix sets form a
set of measure zero within the set of reducible matrices.
\end{conjecture}

In \cite{JSREC} we delve further this
analysis especially explaining how reducible families can be
handled.

\ifpdf
\begin{figure}
\centering
\includegraphics[width=0.6\textwidth]{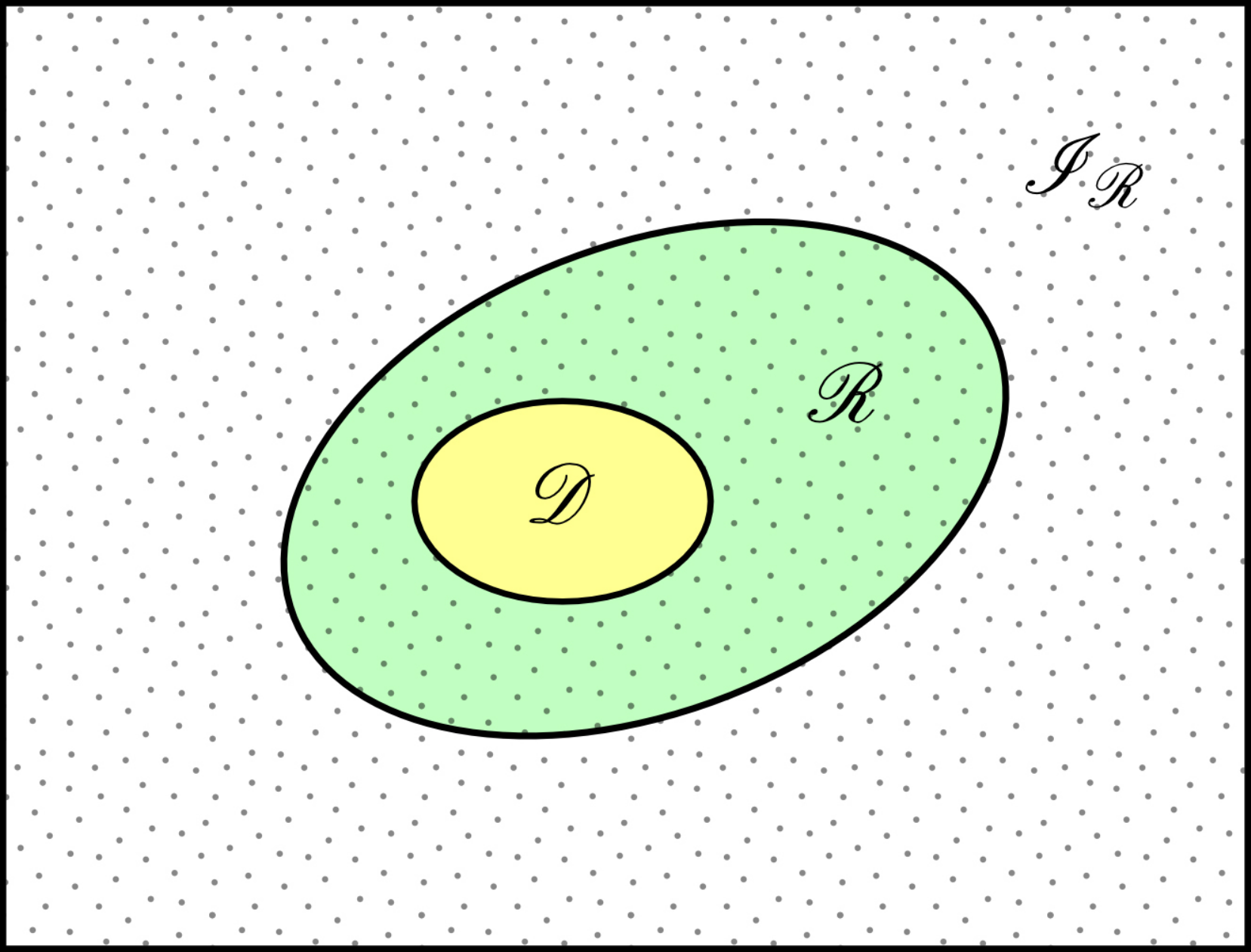}
\caption{Space of bounded families of matrices ${\cal B}$: the set
of defective families is denoted by ${\cal D}$, while its
complement, highlighted by dots, is the set of nondefective families
${\cal N_D}$.} \label{fig:Mspace}
\end{figure}
\else \fi

We add just that Brayton and Tong in \cite{BraTon} give an
alternative sufficient--condition for nondefectiveness. They prove
that, considered each matrix $P$ in the semigroup $\PF$ and the
associated similarity matrix $S_P$ that reduce $P$ into its Jordan
form, if every $S_P$ has columns linearly independent uniformly on
all $P\in\PF$, then $\F$ is nondefective. This alternative
sufficient--condition represents the generalization of the concept
of strongly nondefectiveness to the case of sets of matrices, in
fact for a single matrix $A$ strongly nondefectiveness is equivalent
to semisimplicity of all the eigenvalues in the spectrum of $A$ or
equivalently to diagonalizability of $A$ (ref pages
\pageref{defective1} and \pageref{defective2}). Clearly strongly
nondefectiveness implies nondefectiveness, but checking this
sufficient--condition is not feasible in practice.

As previously mentioned there are not known algorithms for deciding
uniform asymptotic stability of a generic set of matrices and it is
unknown if this problem is algorithmically decidable in general. We
have also seen that uniform asymptotic stability of the set $\F$ is
equivalent to $\rhoF<1$. In order to check if $\rhoF<1$ for finite
families we may think of using the four members inequality
(\ref{eq:fourMember})
\begin{equation*}
    \overline\rho_k(\F)\leq \rho(\F) \leq \widehat\rho_k(\F)\ \ \ \textrm{ for all }  k \geq 1
\end{equation*}
The procedure could be the following \cite{DauLag1}:

\begin{itemize}

\item[(1) ] We evaluate
\begin{equation*}
\overline\rho_k(\F):=\max_{P\in \P_k(\F)} \rho(P)^{1/k} \quad{\rm
and }\quad \widehat\rho_k(\F):=\max_{P\in \P_k(\F)} \|P\|^{1/k}
\end{equation*}
for increasing values of $k\geq 1$.

\item[(2) ] As soon as $\widehat\rho_k< 1$ or $\overline\rho_k\geq 1$ we stop
and declare the set uniform asymptotic stable or unstable
respectively. \end{itemize}

We observe that this procedure always stops after finitely many
steps unless $\rho= 1$ and $\overline\rho_k < 1$ for all $k\geq 1$,
but this never occurs for families that, satisfying the finiteness
property, have an s.m.p.

\begin{definition}[Finiteness property and s.m.p.]\label{def:finitness&smp}
A finite family $\F$ of $n\times n$--matrices has the
\emph{finiteness property} \index{finiteness property} if there
exists, for some $k \geq 1$, a product $\overline{P} \in {\mathcal
P}_k\left(\F \right)$ such that
\begin{equation*}
\rho\left( \overline{P} \right) = \rho\left(\F \right)^k.
\end{equation*}

The product $\overline{P}$ is said to be a
\index{spectrum--maximizing product}\emph{spectrum--maximizing
product} or \emph{s.m.p.}\index{s.m.p.|see{spectrum--maximizing
product}} for~$\F$. Some authors refer to \emph{optimal product}
instead of s.m.p.\index{optimal product}, see for instance
\cite{JuBl08,Mae3}.

An s.m.p. is said \index{s.m.p.!minimal}{\em minimal} if it is not a
power of another s.m.p. of $\F$.

Any eigenvector $x \ne 0$ of $\overline{P}$ related to an eigenvalue
$\lambda$ with $|\lambda|=\rho(\overline{P})$ is said to be a
\index{leading eigenvector}{\em leading eigenvector} of $\F$.
\end{definition}

From the previous definition is evident that uniform asymptotic
stability is algorithmically decidable for finite sets of matrices
that have the finiteness property.

Lagarias and Wang  in $1995$ \cite{LagWan} conjectured that the
finiteness property was valid for all finite families of real
matrices (the so--called \emph{finiteness
conjecture}\index{finiteness conjecture}). Unfortunately this
conjecture does not hold true: Bousch and Mairesse \cite{BouMer} and
later other authors \cite{Blondel1,Koz05} presented
non--constructive counterproofs. In particular in \cite{Blondel1}
Blondel et al. proved that for the parametric family
\begin{equation*}
\mathcal{F}_\alpha=\left\{A,\; \alpha B\right\}=\left\{\left[
  \begin{array}{cc}
    1& 1  \\
   0& 1  \\
  \end{array}
\right],\; \alpha\left[
  \begin{array}{cc}
    1& 0  \\
   1& 1  \\
  \end{array}
\right]\right\}\quad \textrm{with}\quad \alpha \in [0, 1]
\end{equation*}
there exist uncountably many values of the parameter $\alpha$ for
which $\mathcal{F}_\alpha$ does not satisfy the finiteness
conjecture. They were unable to find a single explicit value of
$\alpha$ and they conjectured that the set of values $\alpha \in [0,
1]$ for which the finiteness conjecture is not satisfied is of
measure zero. Recently Hare et al. \cite{MorThe}, using
combinatorial ideas and ergodic theory, have been able to
approximate, up to a desired precision, an explicit value $\alpha$
such that $\mathcal{F}_\alpha$ does not satisfy the finiteness
conjecture. The question if there exist families of matrices with
rational entries that violate the conjecture remains still open.
Based on all the numerical experiments developed in the last years
and the results previously mentioned a new conjecture has been
introduced:
\begin{conjecture}[{\cite[Blondel et al. and
Maesumi]{Blondel1,Mae3}}]\label{conj:finitPropMeasure} The
finiteness property is true a.e. in the space of finite families of
complex square matrices, i.e. the set of families of matrices for
which the finiteness property is not true has measure zero in the
space of finite families.
\end{conjecture}
If this conjecture is true then it suggests us to track s.m.p.'s
candidates out and validate them with some procedure in order to
find the spectral radius of the family. In \cite{JSREC} we
explain how to perform the validation step using particular extremal
norms for the given set.

The idea behind this last conjecture is that the NP--hardness,
non--algebraicity and undecidability results are due to certain rare
and extreme cases and that in the generic case the evaluation of the
spectral radius, while could be computationally intensive, is
possible. About the computational complexity we remind an example,
given by Berger and Wang \cite[Example 2.1]{BerWan}, of a set of two
$\2x2$--matrices with minimal s.m.p. of length $k\geq1$ with $k$
arbitrarily large:

\begin{equation*}
\F=\left\{\alpha^k\left(
     \begin{array}{cc}
       0 & 0 \\
       1 & 0 \\
     \end{array}
   \right), \alpha^{-1}\left(
     \begin{array}{cc}
        \cos\frac{\pi}{2k} & \sin\frac{\pi}{2k} \\
       -\sin\frac{\pi}{2k} & \cos\frac{\pi}{2k} \\
     \end{array}
   \right)\right\}
   \quad \textrm{with} \quad 1<\alpha<\left(\cos\frac{\pi}{2k}\right)^{-1}
\end{equation*}
They prove that $\rhoF=1$, $\overline\rho_j(\F)<1$ for $j\leq k$ and
$\overline\rho_{k+1}(\F)=1$.

We recall that Blondel and Tsitsiklis in \cite{BlTs00b} proved also
that the \emph{effective finiteness conjecture} is false:
\begin{conjecture}[Effective finiteness conjecture]\label{EFC}
For any finite set $\F$ of square matrices with rational entries
there exists an effectively computable natural number $t(\F)$ such
that $\overline\rho_{t(\F)}(\F)=\rhoF$
\end{conjecture}
The falseness of this conjecture implies that, given a family of
matrices with rational entries which admits a spectrum--maximizing
product, the length of the s.m.p. can be arbitrary long and
consequently the computation of the spectral radius can become a
tough problem. Nevertheless for random families this product appears
to be, luckily, quite short in general.

The finiteness property is known to hold in many cases:
\begin{itemize}
  \item when the matrices in $\F$ are all simultaneously
upper--triangularizable, or they can be simultaneously normalized,
or the Lie algebra associated with the set $\F$ is solvable. In
these cases, in fact, the spectral radius is simply equal to the
maximum of the spectral radii of the matrices (Property 10, special
case \ref{item:maxRho} on page \pageref{item:maxRho});
  \item when a finite set of real matrices admits an extremal \emph{piecewise analytic norm} in
$\R^n$.\index{piecewise analytic norm} A piecewise analytic norm is
any norm on $\R^n$ whose unit ball $B$ has a boundary which is
contained in the zero set of a holomorphic function $f(z)$, i.e.
complex differentiable at every point in its domain, defined on a
connected open set $\Omega\in\Cn$ containing $0$, which has
$f(0)\neq 0$ (Lagarias and Wang \cite{LagWan});
  \item when a finite set of real matrices admits an extremal \emph{piecewise algebraic norm}
in $\R^n$.\index{piecewise algebraic norm} A piecewise algebraic
norm is one whose boundary is contained in the zero set of a
polynomial $p(z) \in \R[z_l, \ldots , z_n]$, which has $p(0) \neq
0$. This is the case when the unit ball of a norm is a polytope \cite{JSREC}, or an ellipsoid (ref page \pageref{ellipsNorm}),
or the $l^p$ norm for rational $p$, with $1 \leq p \leq \infty$
(Lagarias and Wang in \cite{LagWan} extended the result proved by
Gurvits in \cite{Gurv} for real polytope extremal norms to the
general case of piecewise algebraic norms in $\R^n$);
  \item when a finite set of matrices admits a complex polytope extremal
norm. This it has been proved by Guglielmi, Wirth and Zennaro in
\cite[Theorem 5.1]{GWZ05} extending to the complex case the results
by Gurvits \cite{Gurv} and Lagarias and Wang \cite{LagWan}. We come
back to polytope norms in \cite{JSREC}.\label{item:bcpFiniteness}
\end{itemize}

For other classes of sets of matrices the finiteness property has
been only conjectured to be true, an example is the class of sets of
matrices with rational entries. Indeed the proof of the finiteness
property for sets of rational matrices would be satisfactory for
practical applications: the matrices that one handles or enters in a
computer are rational--valued.

Recently Blondel and Jungers \cite{JuBl08} have proved the following
Theorem:

\begin{theorem}[{\cite[Theorem 4]{JuBl08}}]\label{th:ratVSsign} $\,$
\begin{enumerate}
  \item The finiteness property holds for all sets of nonnegative rational
matrices if and only if it holds for all pairs of binary matrices.
  \item The finiteness property holds for all sets of rational matrices if
and only if it holds for all pairs of matrices with entries in $\{
-1,0,+1 \}$.
\end{enumerate}
\end{theorem}

They proposed, consequently, the following conjecture
\begin{conjecture}[{\cite[Blondel, Jungers and Protasov]{BJP06,JuBl08}}]
Pairs of binary matrices have the finiteness property.
\end{conjecture}

If this conjecture is correct then, by Theorem \ref{th:ratVSsign},
nonnegative rational matrices also satisfy the finiteness property
and, thus, the question $\rhoF<1$ becomes decidable for sets of
matrices with nonnegative rational entries. From a decidability
point of view this last result would be somewhat surprising since it
is known that the closely related question $\rhoF\leq 1$ is known to
be no algorithmically decidable for such sets of matrices (ref
 Property 13 on page \pageref{item:undec}). Blondel and
Jungers \cite{JuBl08} proved that pairs of $\2x2$ binary--matrices
satisfy the finiteness property and observed that the length of the
s.m.p.'s is always very short. This result is promising even though
a generalization to the case of $\nxn$--matrices seems quite
difficult due to the falseness of the effective finiteness
conjecture \ref{EFC}, which implies that the length of the s.m.p.'s
for families of $\nxn$--matrices can become extremely long.

A more general version of the previous Conjecture is the following
\begin{conjecture}[{\cite[Blondel, Jungers and
Protasov]{BJP06,JuBl08}}]\label{conj:BJP} The finiteness property
holds for pairs of matrices with entries in $\{ -1,0,+1 \}$ (the
so--called sign--matrices).
\end{conjecture}
This last would imply, by Theorem \ref{th:ratVSsign}, that the
finiteness property holds for all sets of rational matrices. In \cite{JSREC} we prove analytically the finiteness property for
pairs of $\2x2$ sign--matrices, i.e. matrices with entries in $\{
-1,0,+1 \}$.


\begin{thebibliography}{19}

\bibitem{Bara88}
N.~E.~Barabanov. Lyapunov indicator of discrete inclusions.
{I--III}. {\it Autom.~Remote Control}, vol. 49 (2,3,5), 1988, pp
152--157, 283--287, 558--565.

%\bibitem{BCS78}
%R.H.~Bartels, A.R~Conn, and J.W.~Sinclair.
%{Minimization techniques for piecewise differentiable functions: the
%$l\sb{1}$ solution to an overdetermined linear system}.
%{\it SIAM J. Numer. Anal.}, vol. 15(2), 1978, pp 224--241.

\bibitem{BerWan}
M.~A.~Berger, Y.~Wang. {Bounded semigroups of matrices}. {\it Linear
Algebra and its Applications}, vol. 166, 1992, pp 21--27.

\bibitem{Blondel}
V.~D.~Blondel. {The birth of the joint spectral radius: An interview
with Gilbert Strang}. {\it Linear Algebra and its Applications},
vol. 428, 2008, pp 2261--2264.

\bibitem{BloCan03}
V.~D.~Blondel, V.~Canterini. {Undecidable problems for probabilistic
automata of fixed dimension}. {\it Theory of Computing Systems},
vol. 36 (3), 2003, pp 231--245.

\bibitem{BloChan10}
V.~D.~Blondel, C.~T.~Chang. {A comparison of extremal norms methods
to approximate the joint spectral radius}. {\it Submitted to 49th
IEEE Conference on Decision and Control (CDC10)}, 2010.

\bibitem{BJP06}
V.~D.~Blondel, R.~Jungers, V.~Yu.~Protasov. {On the complexity of
computing the capacity of codes that avoid forbidden difference
patterns}. {\it IEEE Trans. Inform. Theory}, vol. 52, 2006, pp
5122--5127.

\bibitem{BJP10}
V.~D.~Blondel, R.~Jungers, V.~Yu.~Protasov. {Joint spectral
characteristics of matrices: a conic programming approach}. {\it
SIAM J. Matrix Anal. Appl.}, vol. 31 (4), 2010, pp 2146--2162.

\bibitem{Blondel2}
V.~D.~Blondel, Y.~Nesterov. {Computationally efficient
approximations of the joint spectral radius}. {\it SIAM J. Matrix
Analysis Appl.}, vol. 27 (1), 2005, pp 256--272.

\bibitem{Blondel3}
V.~D.~Blondel, Y.~Nesterov, J.~Theys. {On the accuracy of the
ellipsoid norm approximation of the joint spectral radius}. {\it
Linear Algebra and its Applications}, vol. 394 (1), 2005, pp
91--107.

\bibitem{Blondel1}
V.~D.~Blondel, J.~Theys, A.~A.~Vladimirov. {An elementary
counterexample to the finiteness conjecture}. {\it SIAM J. Matrix
Analysis Appl.}, vol. 24 (4), 2003, pp 963--970.

\bibitem{BlTs99}
V.~D.~Blondel, J.~N.~Tsitsiklis. Complexity of stability and
controllability of elementary hybrid systems. {\it Automatica}, vol.
35 (3), 1999, pp 479--489.

\bibitem{BlTs00}
V.~D.~Blondel, J.~N.~Tsitsiklis. A survey of computational
complexity results in systems and control. {\it Automatica}, vol. 36
(9), 2000, pp 1249--1274.

\bibitem{BlTs00b}
V.~D.~Blondel, J.~N.~Tsitsiklis. The boundedness of all products of
a pair of matrices is undecidable. {\it Systems \& Control Letters},
vol. 41 (2), 2000, pp~135--140.

\bibitem{BouMer}
T.~Bousch, J.~Mairesse. {Asymptotic height optimization for topical
IFS, Tetris heaps and the finiteness conjecture}. {\it J. Amer.
Math. Soc.}, vol. 15, 2002, pp 77--111.

\bibitem{BraTon}
B.~K.~Brayton, C.~H.~Tong. {Constructive stability and asymptotic
stability of dynamical systems}. {\it IEEE Trans. Circuits Systems},
vol. 27 (11), 1980, pp 1121--1130.

\bibitem{CavaMi}
A. S. Cavaretta, C. A. Micchelli. {The design of curves and surfaces
by subdivision algorithms}. {\it Mathematical methods in computer
aided geometric design} Academic Press, Boston, MA, 1989, pp
115--153.

\bibitem{CheZho}
Q.~Chen, X.~Zhou. {Characterization of joint spectral radius via
trace}. {\it Linear Algebra and its Applications}, vol. 315, 2000,
pp 175--188.

\bibitem{JSREC}
A. Cicone. On the Joint Spectral Radius Exact Computation.

\bibitem{artic2}
A. Cicone, N. Guglielmi, S. Serra--Capizzano, M. Zennaro. Finiteness
property of pairs of $\2x2$ sign--matrices via real extremal
polytope norms. {\it Linear Algebra and its Applications}, vol. 432
(2--3), 2010, pp 796--816.

\bibitem{CoHe}
D. Colella, C. Heil. {The characterization of continuous,
four--coefficient scaling functions and wavelets}. {\it IEEE Trans.
Inform. Theory}, vol. 38 (2, Part 2), 1992, pp 876--881.

\bibitem{DauLag1}
I.~Daubechies, J.~C.~Lagarias. {Sets of matrices all infinite
products of which converge}. {\it Linear Algebra and its
Applications}, vol. 161, 1992, pp 227--263.

\bibitem{DauLag2}
I.~Daubechies, J.~C.~Lagarias. {Corrigendum/addendum to: Sets of
matrices all infinite products of which converge}. {\it Linear
Algebra and its Applications}, vol. 327, 2001, pp 69--83.

\bibitem{Els}
L.~Elsner. {The generalized spectral--radius theorem: an
analytic--geometric proof}. {\it Linear Algebra and its
Applications}, vol. 220, 1995, pp 151--159.

\bibitem{GhoshLee}
S. Ghosh, J. W. Lee. Equivalent conditions for uniform asymptotic
consensus among distributed agents. {\it American Control
Conference}, 2010, pp 4821--4826.

\bibitem{Gri96}
G.~Gripenberg. {Computing the joint spectral radius}. {\it Linear
Algebra and its Applications}, vol. 234, 1996, pp 43--60.

\bibitem{GMV10}
N.~Guglielmi, C.~Manni, D.~Vitale. Convergence analysis of $C^2$
Hermite interpolatory subdivision schemes by explicit joint spectral
radius formulas. {\it Linear Algebra and its Applications}, vol. 434
(4), 2011, pp 884--902.

\bibitem{GWZ05}
N.~Guglielmi, F.~Wirth, M.~Zennaro. {Complex polytope extremality
results for families of matrices}. {\it SIAM J. Matrix Analysis
Appl.}, vol. 27, 2005, pp 721--743.

\bibitem{GugZen00}
N.~Guglielmi, M.~Zennaro. {On the asymptotic properties of a family
of matrices}. {\it Linear Algebra and its Applications}, vol. 322,
2001, pp 169--192.

\bibitem{GugZen01}
N.~Guglielmi, M.~Zennaro. {On the zero--stability of variable
stepsize multistep methods: the spectral radius approach}. {\it
Numer. Math.}, vol. 88, 2001, pp 445--458.

\bibitem{GugZen03}
N.~Guglielmi, M.~Zennaro. {On the limit products of a family of
matrices}. {\it Linear Algebra and its Applications}, vol. 362,
2003, pp 11--27.

\bibitem{GugZen03b}
N.~Guglielmi, M.~Zennaro. {Stability of one--leg $\Theta$--methods
for the variable coefficient pantograph equation on the
quasi--geometric mesh}. {\it IMA J. Numer. Anal.}, vol. 23, 2003, pp
421--438.

\bibitem{GugZen05b}
N.~Guglielmi, M.~Zennaro. {Polytope norms and related algorithms for
the computation of the joint spectral radius}. {\it Proceedings of
the 44th IEEE Conference on Decision and Control and European
Control Conference CDC--ECC'05}, 2005, pp 3007--3012.

\bibitem{GugZen05}
N.~Guglielmi, M.~Zennaro. {Balanced complex polytopes and related
vector and matrix norms}. {\it J. Convex Anal.}, vol. 14, 2007, pp
729--766.

\bibitem{GugZen08}
N.~Guglielmi, M.~Zennaro. {An algorithm for finding extremal
polytope norms of matrix families}. {\it Linear Algebra and its
Applications}, vol. 428, 2008, pp 2265--2282.

\bibitem{Gurv}
L.~Gurvits. {Stability of discrete linear inclusions}. {\it Linear
Algebra and its Applications}, vol. 231, 1995, pp 47--85.

\bibitem{Gurv2}
L.~Gurvits. Stability of Linear Inclusions--Part 2, {\it NECI
Technical Report TR}, 1996, pp 96--173.

\bibitem{MorThe}
K.~G.~Hare, I.~D.~Morris, N.~Sidorov, J.~Theys. {An explicit
counterexample to the Lagarias--Wang finiteness conjecture}. {\it To
appear}. {\url{http://arxiv.org/abs/1006.2117}}

\bibitem{HeiStra}
C.~Heil, G.~Strang. {Continuity of the joint spectral radius:
application to wavelets}. {\it Linear Algebra for Signal Processing
, A. Bojanczyk and G. Cybenko (eds.), IMA Vol. Math. Appl.,
Springer, New York}, vol. 69, 1995, pp 51--61.

\bibitem{Heu}
H. G.~Heuser. {\it Functional analysis}. John Wiley \& Sons, New
York, 1982.

\bibitem{horn}
R. A. Horn, C. R. Johnson. {\it Matrix Analysis}. Cambridge
University Press, 2007.

\bibitem{JadLinMor03}
A. Jadbabaie, J. Lin, A. S. Morse. Coordination of groups of mobile
autonomous agents using nearest neighbor rules. {\it IEEE Trans.
Automat. Control}, vol. 48 (6), 2003, pp 988--1001.

\bibitem{Ju09}
R.~Jungers. {\it The joint spectral radius, theory and
applications}. In Lecture Notes in Control and Information Sciences,
Springer--Verlag, Berlin, vol. 385, 2009.

\bibitem{JuBl08}
R.~Jungers, V.~Blondel. {On the finiteness properties for rational
matrices}. {\it Linear Algebra and its Applications}, vol. 428,
2008, pp 2283--2295.


\bibitem{JuPr09}
R.~Jungers, V.~Yu.~Protasov. {Counterexamples to the CPE
conjecture}. {\it SIAM J. Matrix Analysis Appl.}, vol. 31, 2009, pp
404--409.

\bibitem{Koz90}
V. S.~Kozyakin. {Algebraic unsolvability of problem of absolute
stability of desynchronized systems}. {\it Automat. Remote Control},
vol. 51 (6), 1990, pp 754--759.

\bibitem{Koz05}
V. S.~Kozyakin. {A dynamical systems construction of a
counterexample to the finiteness conjecture.}. {\it Proceedings of
the 44th IEEE Conference on Decision and Control and European
Control Conference CDC--ECC'05}, 2005, pp 2338--2343.

\bibitem{Koz10}
V. S.~Kozyakin. {An explicit Lipschitz constant for the joint
spectral radius}. {\it Linear Algebra and its Applications}, vol.
433, 2010, pp 12--18.

\bibitem{LagWan}
J. C.~Lagarias, Y.~Wang. {The finiteness conjecture for the
generalized spectral radius of a set of matrices}. {\it Linear
Algebra and its Applications}, vol. 214, 1995, pp 17--42.

\bibitem{Lev}
J.~Levitzki. {Über nilpotente Unterringe}. {\it Math. Ann.}, vol.
105, 1931, pp 620--627.

\bibitem{Mae1}
M.~Maesumi. {Optimum unit ball for joint spectral radius: an example
from  four--coefficient MRA}. {\it in Approximation Theory VIII:
Wavelets and Multilevel Approximation. C.K. Chui and L.L. Schumaker
(eds.)},  vol. 2, 1995, pp 267--274.

\bibitem{Mae2}
M.~Maesumi. {Calculating joint spectral radius of matrices and
H\"{o}lder exponent of  wavelets}. {\it in Approximation Theory IX.
C.K. Chui and L.L. Schumaker  (eds.)}, 1998, pp 1--8.

%\bibitem{Mae3}
%M.~Maesumi. {Construction of optimal norms for semigroups of matrices}.
%{\it Proceedings of the 44th IEEE Conference on Decision and Control and
%European Control Conference CDC--ECC'05}, 2005, pp 3013--3018.

\bibitem{Mae3}
M.~Maesumi. {Optimal norms and the computation of joint spectral
radius of matrices}. {\it Linear Algebra and its Applications}, vol.
428 (10), 2008, pp 2324--2338.

\bibitem{MiPra89}
C. A. Micchelli, H. Prautzsch. {Uniform refinement of curves}. {\it
Linear Algebra and its Applications}, vol.  114--115, 1989, pp
841--870.

\bibitem{MOrSi01}
B. E. Moision, A. Orlitsky, P. H. Siegel. On codes that avoid
specified differences. {\it IEEE Trans. Inform. Theory}, vol. 47
(1), 2001, pp 433--442.

\bibitem{Moreau}
L. Moreau. Stability of multi--agent systems with time--dependent
communication links. {\it IEEE Trans. Automat. Control}, vol. 50
(2), 2005, pp 169--182.

\bibitem{Moss}
B. M\"ossner. {On the joint spectral radius of matrices of order 2
with equal spectral radius}. {\it Advances in Comp. Math.}, vol. 33
(2), 2010, pp 243--254.

\bibitem{Pat}
M. S. Paterson. {Unsolvability in $3\! \times \!2\,$ matrices}. {\it
Studies in Applied Mathematics}, vol. 49, 1970, pp 105--107.

\bibitem{Pro2}
V.~Yu.~Protasov. {The joint spectral radius and invariant sets of
linear operators}. {\it Fundamentalnaya i prikladnaya matematika},
vol. 2 (1), 1996, pp 205--231.

\bibitem{Pro3}
V.~Yu.~Protasov. {A generalized joint spectral radius, a geometric
approach}. {\it Izvestiya: Mathematics}, vol. 61 (5), 1997, pp
995--1030.

\bibitem{Pro1}
V.~Yu.~Protasov. {The geometric approach for computing the joint
spectral radius}. {\it Proceedings of the 44th IEEE Conference on
Decision and Control and European Control Conference CDC--ECC'05},
2005, pp 3001--3006.

\bibitem{Pro4}
V.~Yu.~Protasov. {Applications of the Joint Spectral Radius to Some
Problems of Functional Analysis, Probability and Combinatorics}.
{\it Proceedings of the 44th IEEE Conference on Decision and Control
and European Control Conference CDC--ECC'05}, 2005, pp 3025--3030.

%\bibitem{RadRos00}
%R.~Radjavi, P.~Rosenthal. {\it Simultaneous triangularization},
%Universitext, Springer--Verlag, New--York, 2000.

\bibitem{RotStr}
G.~C.~Rota, G.~Strang. {A note on the joint spectral radius}. {\it
Indagatione Mathematicae}, vol. 22, 1960, pp 379--381.
{\url{http://books.google.it/books?id=x7zKLGu3g9IC&lpg=RA1-
PA40&ots= rxuCbvNTxI&dq=a%20note%20on%20the%20joint%20spectral%20
radius%20rota%20strang&lr&pg=PA74#v=onepage&q&f=false}}

\bibitem{Shi}
M.~H.~Shih. {Simultaneous Schur stability}. {\it Linear Algebra and
its Applications}, vol 287, 1999, pp 323--336.

\bibitem{ShWuPa}
M.~H.~Shih, J.~W.~Wu, C.~T.~Pang. {Asymptotic stability and
generalized Gelfand spectral radius formula}. {\it Linear Algebra
and its Applications}, vol. 252, 1997, pp 61--70.

\bibitem{Str}
G.~Strang. The joint spectral radius, Commentary on paper n. 5. {\it
Collected Works of Gian--Carlo Rota}, 2001.
{\url{http://www-math.mit.edu/~gs/papers/ss.pdf}}

\bibitem{Theys}
J.~Theys. {Joint Spectral Radius : Theory and approximations}. {\it
PhD thesis}. Université catholique de Louvain, 2005.
{\url{http://www.inma.ucl.ac.be/~blondel/05thesetheys.pdf}}

\bibitem{Tsi87}
J.~N.~Tsitsiklis. On the stability of asynchronous iterative
processes. {\it Math. Systems Theory }, vol. 20 (2--3), 1987, pp
137--153.

\bibitem{TsitBlon97}
J.~N.~Tsitsiklis, V.~D.~Blondel. The Lyapunov exponent and joint
spectral radius of pairs of matrices are hard -- when not impossible
-- to compute and to approximate. {\it Mathematics of Control,
Signals, and Systems}, vol. 10, 1997, pp 31--40.

\bibitem{vdb}
R.~J.~Vanderbei. {\it Linear programming, Foundations and
extensions}. International Series in Operations Research \&
Management Science, 114. Third edition. Springer, New York, 2008.

\bibitem{Wirt02}
F.~Wirth. The generalized spectral radius and extremal norms. {\it
Linear Algebra and its Applications}, vol. 342, 2002, pp 17--40.

\bibitem{Wirt05}
F.~Wirth. The generalized spectral radius is strictly increasing.
{\it Linear Algebra and its Applications}, vol. 395, 2005, pp
141--153.

\bibitem{Zie}
G. M.~Ziegler. {\it {Lectures on polytopes}}. Springer--Verlag, New
York, 1995.

\end{thebibliography}
\end{document}